
\input amstex
\documentstyle{amsppt}
\magnification=\magstep1
\def\bc{{\Bbb C}} 
\def\bp{{\Bbb P}}

\def\dd{{D}}  

\def\sing{\operatorname{Sing}}
\def\hip{{\bc\bp^n_\infty}}
\def\hi{{\bc\bp^2_\infty}}

\newcount\parno \parno=1
\newcount\prono \prono=1
\def\sec{\S\the\parno.-\ \global\prono=1}
\def\etiqueta{\hbox{(\the\parno.\the\prono)}}
\def\finparrafo{\global\advance\parno by1
\vskip.1truecm\ignorespaces}
\def\finparrafo{\global\advance\parno by1
\vskip.1truecm\ignorespaces}
\def\cita{\ignorespaces\ \the\parno.\the\prono%
\global\advance\prono by 1}

\NoBlackBoxes

\document

\NoBlackBoxes
\topmatter

\title
On the zeta-function of a polynomial at infinity
\endtitle

\author
S.M. Gusein--Zade, I. Luengo, A. Melle--Hern\'andez
\endauthor

\address
Moscow State University, 
Faculty of Mathematics and Mechanics,
Moscow, 119899, Russia.
\endaddress
\email
sabir\@ium.ips.ras.ru
\endemail
\address
Departamento de \'Algebra,
Universidad Complutense,
Ciudad Universitaria s/n,
E-28040 Madrid, Spain.
\endaddress
\email
iluengo\@eucmos.sim.ucm.es
\endemail
\address
Departamento de Geometr\'{\i}a y Topolog\'{\i}a,
Universidad Complutense,
Ciudad Universitaria s/n,
E-28040 Madrid, Spain.
\endaddress
\email
amelle\@eucmos.sim.ucm.es
\endemail

\thanks
First author was partially supported by Iberdrola, INTAS--96--0713, RFBR 96--15--96043.
Last two authors were partially supported by DGCYT PB94-0291.
\endthanks

\abstract
We use the notion of Milnor fibres of the germ of a meromorphic function
and the method of partial resolutions for a study of topology of a polynomial 
map at infinity  (mainly for calculation of the zeta-function of a monodromy).
It gives effective  methods of computation of the zeta-function for a number 
of cases and a criterium for a value to be atypical at infinity. 
\endabstract
\endtopmatter

\document

\head\sec Introduction
\endhead

The main idea of the paper is to bring together
methods of \cite{7} and \cite{8} for computing the 
zeta-function of the monodromy at infinity of a polynomial.
Let $P$ be a complex polynomial in $(n+1)$ variables. It defines a map from
$\bc^{n+1}$ to $\bc$ which also will be denoted by $P.$ It is known (\cite{13})
that there exists a finite set $B(P)\subset \bc$ such that
the map $P$ is a $C^\infty$ locally trivial fibration over its
complement. 
The monodromy
transformation $h$ of this fibration corresponding to the loop $z_0\cdot\exp(2\pi i
\tau)$ ($0\leq\tau\leq 1$) with $\|z_0\|$ big enough is called the {\it geometric monodromy 
at infinity} of the polynomial $P.$ 
Let $h_{*}$ be its action in the homology groups of the fibre (the
level set) $\{P=z_0\}$.

\definition{Definition} The {\it zeta-function of the monodromy at infinity} of the
polynomial $P$ is the rational function
$$\zeta_P(t)=
\prod_{q\geq 0} \{ \det\,[\,id-t\,h_{*}|_{H_q(\{P=z_0\};\bc)}]\}^{(-1)^q}.$$
\enddefinition

\remark{Remark 1} We use the definition from \cite{2}, which means
that the zeta-function defined this way is the inverse of that used in \cite{1}.
\endremark

The degree of the zeta-function (the degree of the numerator minus the degree of
the denominator) is equal to the Euler characteristic $\chi_P$ of the (generic)
fibre $\{P=z_0\}.$ Formulae for the zeta-functions at infinity for certain
polynomials were given in particular in \cite{6}, \cite{9}. 

\bigbreak
\finparrafo

\head \sec Zeta-function of a polynomial via 
zeta-functions of meromorphic germs
\endhead

A polynomial function $P:\bc^{n+1}\to\bc$ defines a meromorphic function $P$
 on the projective space $\bc\bp^{n+1}.$ At each point $x$ of the infinite hyperplane
$\bc\bp^n_\infty$ the germ of the meromorphic function $P$ has
the form
$\displaystyle\frac{F(u,x_1,\ldots,x_n)}{u^d}$ where $u,x_1,\ldots,x_n$ are local
coordinates such that $\hip=\{u=0\},$ $F$ is the germ of a holomorphic function, and $d$ is the
degree of the polynomial
$P.$

In \cite{8}, for a meromorphic germ  $f=\frac{F}{G},$ there were defined two
Milnor fibres (the zero and the infinite ones), two monodromy transformations,
and thus two zeta-functions $\zeta_f^0(t)$ and $\zeta_f^\infty(t)$.
Let $\zeta^{\bullet}_{P,x}(t)$ ($\bullet=0$ or $\infty$) be the corresponding
zeta-function of the germ of the meromorphic function $P$ at the point $x\in\hip.$

For the aim of convinience, in \cite{8} we considered only meromorphic germs 
$f=\frac{F}{G}$ with $F(0)=G(0)=0.$ At a generic point of the infinite hyperplane $\hip$ the
meromorphic function $P$ has the form $\frac{1}{u^d}.$ For a germ of the form
$f=\frac{1}{G}$ with $G(0)=0,$ it is reasonable to give the following definition: its infinite
Milnor fibre coincides with the (usual) Milnor fibre of the holomorphic germ $G$
and its zero Milnor fibre is empty. Thus $\zeta_f^0(t)=1$ and $\zeta_f^\infty(t)=
\zeta_G(t).$ According to this definition, for the germ $\frac{1}{u^d},$ its
infinite zeta-function is equal to
$(1-t^d).$

Let ${\Cal S}=\{ \Xi\}$ be a
prestratification of the infinite hyperplane $\hip$ (that is a partitioning of $\hip$ into
semi-analytic subspaces without any regularity conditions) such that,
for each stratum $\Xi$ of $\Cal S,$ the infinite zeta-function
$\zeta_{P,x}^\infty(t)$ does not depend on $x,$ for $x \in \Xi$. Let us
denote this zeta-function by $\zeta_{\Xi}^\infty(t)$ and by
$\chi_{\Xi}^\infty$  its degree
$\deg\zeta_{\Xi}^\infty(t)$. A straighforward repetition of the
arguments from the proof of Theorem 1 in \cite{7} gives

\proclaim{Theorem 1} 
$$\zeta_{P}(t)=\prod_{\Xi\in {\Cal S}}[\zeta_{\Xi}^\infty(t)]^{\chi (\Xi)},$$
$$
\chi_{P}=\sum_{\Xi\in {\Cal S}}\chi_{\Xi}^\infty\cdot\chi(\Xi).$$
\endproclaim

\remark{Remark 2} One can write the formula for $\chi_P$ in the form of an 
integral with respect to the
Euler characteristic 
$$\chi_P=\int_\hip \chi_{P,x}^\infty \,d\chi$$
\noindent in the sense of Viro (\cite{14}).
\endremark

\remark{ Remark 3} Let $P_d$ be the (highest) homogeneous part of
degree
$d$ of the polynomial $P.$ Then at each point $x\in \hip\setminus \{P_d=0\}$
the germ of the meromorphic function $P$ is of the form $\frac{1}{u^d}.$
The set $\Xi^n=\hip\setminus \{P_d=0\}$ can be considered as the $n$-dimensional
stratum of a partitioning. It brings the factor $(1-t^d)^{\chi(\Xi^n)}$ into
the zeta-function $\zeta_P(t).$ 
\endremark
\bigbreak
\finparrafo

\head\sec Examples
\endhead
\medbreak
\subhead 3.1. Yomdin-at-infinity polynomials
\endsubhead
This name was introduced in \cite{4}.
For a polynomial $P\in\bc[z_0,z_1,\ldots,z_n],$ let $P_i$ be its homogeneous part of
degree $i.$ Let a polynomial $P$ be of the form $P=P_d+P_{d-k}+$\,terms of lower
degree, $k\geq 1.$ Let us consider
hypersurfaces
in $\bc\bp^n$ defined
by  $\{P_d=0\}$ and $\{P_{d-k}=0\}$.
Let $\sing(P_d)$ be the singular locus of the hypersurface $\{P_d=0\}$
(including all points where $\{P_d=0\}$ is not reduced). One says that  $P$
is a {\it Yomdine-at-infinity polynomial\/} if
$\sing(P_d)\cap \{P_{d-k}=0\}=\emptyset$ (in particular it implies that
$\sing(P_d)$ is finite).

Y. Yomdin (\cite{15}) has considered critical points of holomorphic functions which
are local versions of such polynomials. He gave a formula for their Milnor
numbers. The generic fibre (level set) of a Yomdin-at-infinity polynomial is
homotopy equivalent to the bouquet of $n$-dimensional spheres (\cite{5}).
Its Euler characteristic $\chi_P$ (or rather the (global) Milnor number) has
been determined in \cite{4}. For $k=1,$ the zeta-function of such a polynomial
has been obtained in \cite{6}.

Let $P(z_0,z_1,\ldots,z_n)=P_d+P_{d-k}+\ldots$ be a Yomdin-at-infinity 
polynomial. Let $\sing(P_d)$ consist of $s$ points $Q_1,\ldots,Q_s.$ One has the
following natural stratification of the infinite hyperplane $\hip$:

\roster
\item  the $n$-dimensional stratum
$\Xi^{n}=\hip\setminus \{P_d=0\};$
\item  the
$(n-1)$-dimensional stratum $\Xi^{n-1}=\{P_d=0\}\setminus \{
Q_1,\ldots,Q_s \}$;
\item  the 0-dimensional strata
$\Xi_i^0$ $(i=1,\ldots,s),$ each consisting of one point $Q_i.$
\endroster

 The Euler characteristic of the stratum $\Xi^{n}$ is equal to 
$$\chi(\hip)-\chi(\{P_d=0\})=(n+1)-\chi(n,d)+(-1)^{n-1}\sum_{i=1}^s\mu_i,$$ 
\noindent where
$\chi(n,d)=(n+1)+\frac{(1-d)^{n+1}-1}{d}$ is the Euler characteristic of
a non-singular hypersurface of degree $d$ in the complex projective space $\hip,$
$\mu_i$ is the Milnor number of the germ of the hypersurface $\{P_d=0\}\subset\hip$ at the
point $Q_i.$ 
At each point of the stratum
$\Xi^{n},$ the germ of the meromorphic function $P$ has (in some local coordinates
$u,y_1,\ldots,y_n$) the form $\frac{1}{u^d}$ ($\hip=\{u=0\}$) and its infinite
zeta-function $\zeta_{\Xi^{n}}^\infty(t)$ is equal to $(1-t^d).$

At each point
of the stratum $\Xi^{n-1},$ the germ of the polynomial $P$ has (in some local
coordinates
$u,y_1,\ldots,y_n$) the form $\frac{y_1}{u^d}.$
Its infinite zeta-function $\zeta_{\Xi^{n-1}}^\infty(t)$ is equal to 1 and thus it
does not contribute a factor to the zeta-function of the polynomial $P.$

At a point $Q_i$ ($i=1,\ldots,s$), the germ of the meromorphic function $P$  has the form 
$\varphi(u,y_1,\ldots,y_n)=\displaystyle\frac{g_i(y_1,\ldots,y_n)+u^k}{u^d},$ where
$g_i$ is a local equation of the hypersurface $\{P_d=0\}\subset\hip$ at the point
$Q_i.$ Thus $\mu_i$ is its Milnor number.

To compute the infinite zeta-function $\zeta_\varphi^\infty(t)$ of the
meromorphic germ $\varphi,$ let us consider a resolution
$\pi:({\Cal X},{\Cal D})\to (\bc^n,0)$  of the singularity $g_i,$ i.e., a
proper modification of $(\bc^n,0)$ which is an isomorphism outside the
origin in $\bc^n$ and such that, at each point of the exceptional divisor 
$\Cal D,$ the lifting $g_i\circ \pi$ of the function  $g_i$ to the space $\Cal X$  of the
modification has (in some local coordinates) the form $y_1^{m_1}\cdot\ldots \cdot
y_n^{m_n}\, (m_i\geq 0)$.  

Let us consider the modification 
${\widetilde \pi}=id\times \pi:(\bc_u\times {\Cal X},0\times {\Cal D})\to
(\bc^{n+1},0)=(\bc_u\times \bc^n,0)$ of the space $(\bc^{n+1},0)$ -- the trivial
extension: $(u,x)\mapsto (u,\pi(x)).$
Let ${\widetilde \varphi}=\varphi\circ {\widetilde \pi}$ be the lifting of the
meromorphic function $\varphi$ to the space $\bc_u\times {\Cal X}$ of the
modification $\widetilde \pi.$ Let ${\Cal M}_{{\widetilde \varphi}}^\infty={\widetilde
\pi}^{-1}({\Cal M}_{\varphi}^\infty)$ (${\Cal M}_{\varphi}^\infty$ is the infinite Milnor
fibre of the germ $\varphi$) be the local level set of the meromorphic function ${\widetilde
\varphi}$ (close to the infinite one).  In the natural way one has the (infinite) monodromy
$h^\infty_{\widetilde \varphi}$ acting on
${\Cal M}_{{\widetilde
\varphi}}^\infty$ and its zeta-function $\zeta_{\widetilde \varphi}^\infty(t).$

\proclaim{Theorem 2}
$$\zeta_{\widetilde \varphi}^\infty(t)=(1-t^{d-k})^{\chi({\Cal D})-1}\,
\zeta_\varphi^\infty (t).$$
\endproclaim

\demo{Proof}  The infinite monodromy transformation of the function $\widetilde\varphi$ can be
described in the following way. Let $h_\varphi^\infty:{\Cal M}_{
\varphi}^\infty \to {\Cal M}_{\varphi}^\infty$ be the infinite monodromy
transformation of the germ
$\varphi.$ One can suppose that it preserves the intersection of the Milnor fibre
${\Cal M}_{\varphi}^\infty$ with the line $\bc_u\times\{0\}.$ There it
coincides  with the infinite monodromy transformation of the
restriction $\varphi|_{\bc_u\times \{0\}}=\frac{u^k}{u^d}$ of the germ
$\varphi$ to this line, i.e., with a cyclic permutation of $(d-k)$ points.
The zeta-function of a cyclic permutation
of $(d-k)$ points is equal to $(1-t^{d-k}).$ The projection 
${\widetilde \pi}:{\Cal M}_{{\widetilde \varphi}}^\infty \to{\Cal M}_{\varphi}^\infty$
is an isomorphism  outside ${\Cal M}_{\varphi}^\infty\cap
(\bc_u\times\{0\})$,
the preimage of each point from ${\Cal M}_{\varphi}^\infty\cap
(\bc_u\times\{0\})$ is isomorphic to the exceptional divisor $\Cal D.$
This means that the transformation (the diffeomorphism) 
${h}_{\widetilde \varphi}^\infty:{\Cal M}_{
\widetilde \varphi}^\infty \to {\Cal M}_{\widetilde \varphi}^\infty$
can be constructed in such a way that it preserves 
${\widetilde \pi}^{-1}({\Cal M}_{\varphi}^\infty\cap (\bc_u\times\{0\}))$ and
acts on it by a cyclic permutation of $(d-k)$ copies of ${\Cal D}.$ The
zeta-function of this transformation of $\{(d-k)  \hbox{ points}\}\times
{\Cal D}$ is equal to 
$(1-t^{d-k})^{\chi({\Cal D})}.$ The result
follows from the {\it multiplication property} of the zeta-function
of a transformation (see \cite{2} p. 94). \qed
\enddemo

For ${\bar m}=(m_1,m_2,\ldots,m_n)$ with integer $m_1\geq
m_2\geq\ldots\geq m_n\geq 0,$ let $S_{\bar m}$ be the set of points
of the exceptional divisor $\Cal D$ of the resolution $\pi$ at which the
lifting of the germ
$g_i$ has the form
$y_1^{m_1}\cdot\ldots \cdot y_n^{m_n};$ for $m\ge 1,$ let 
$S_m$ be $S_{(m,0,\ldots,0)}$. By the formula of A'Campo (\cite{1})
$$\zeta_{g_i}(t)=\prod_{m\ge 1}(1-t^m)^{\chi(S_m)}.\eqno(1)$$
\noindent At a point $x\in\{0\}\times S_{\bar m}\subset\{0\}\times
{\Cal D},$ the lifting ${\widetilde \varphi}=\varphi\circ {\widetilde \pi}$ of the
function $\varphi$ has the local form 
$\displaystyle\frac{y_1^{m_1}\cdot\ldots \cdot y_n^{m_n}+u^k}{u^d}.$ Thus, for fixed
$\bar m$, the infinite zeta-function 
$\zeta_{{\widetilde \varphi},x}^\infty(t)$ of the germ of the meromorphic
function $\widetilde
\varphi$ at a point $x$ from $\{0\}\times S_{\bar m}$ is one and the same. It can be
determined by the  Varchenko type formula from \cite{8}.
If there are more than one integers $m_i$ different from zero, 
$\zeta_{{\widetilde \varphi},x}^\infty(t)=(1-t^{d-k}).$ For $x\in \{0\}
\times S_{m},$
$$\zeta_{{\widetilde\varphi},x}^\infty(t)=(1-t^{d-k})
(1-t^{\frac{m(d-k)}{g.c.d.(m,k)}})^{-g.c.d.(m,k)}.$$
\noindent According to Theorem 1
$$\zeta_{{\widetilde\varphi}}^\infty(t)=(1-t^{d-k})
^{\chi({\Cal D})}\prod_{m\geq 1}
\left(1-t^{\frac{m(d-k)}{g.c.d.(m,k)}}\right)^{-g.c.d.(m,k)\cdot\chi(S_m)}$$
\noindent and by Theorem 2
$$\zeta_{\varphi}^\infty(t)=(1-t^{d-k})
\prod_{m\geq 1}\left(1-t^{\frac{m(d-k)}{g.c.d.(m,k)}}\right)^{-g.c.d.(m,k)\cdot\chi(S_m)}.
\eqno(2)$$

The zeta-function $\zeta_h(t)$ of a transformation $h:X\to X$ of a
space
$X$ into itself determines  the zeta-function $\zeta^k_h(t)\,$ of the
$k$-th power $h^k$ of the transformation $h.$   In particular, if $\zeta_h(t)=
\prod\limits_{m\geq 1} (1-t^m)^{a_m},$ then $$\zeta^k_h(t)=\prod_{m\geq
1}
\left(1-t^{\frac{m}{g.c.d.(k,m)}}\right)^{g.c.d.(k,m)\cdot a_m}.$$

The formulae $(1)$ and $(2)$ mean that
$$
\zeta_{\varphi}^\infty (t)=(1-t^{d-k})\left(\zeta_{g_i}^k(t^{d-k})
\right)^{-1}\eqno(3).
$$
\noindent Combining the computations for the stratification 
$\{\Xi^n,\Xi^{n-1},\Xi_i^0\}$ of the infinite hyperplane $\hip$, one has

\proclaim{Theorem 3} For a Yomdin-at-infinity polynomial
$P=P_d+P_{d-k}+\ldots,$ its zeta-function at infinity is equal to
$$\zeta_P(t)=
(1-t^d)^{\chi(\Xi^n)}(1-t^{d-k})^s\left(\prod_{i=1}^s
\zeta_{g_i}^k(t^{d-k})\right)^{-1},$$
\noindent where
$\chi(\Xi^n)=\frac{1-(1-d)^{n+1}}{d}+(-1)^{n-1}\sum_{i=1}^s\mu (g_i)$
and $g_i$ is a local equation of the hypersurface $\{P_d=0\}\subset\hip$
at its singular point $Q_i.$
\endproclaim
\medskip

\subhead 3.2
\endsubhead
Let $(n+1)$ be equal to $3,$ $P=P_d+P_{d-k}+\ldots,\,$ $\{P_d=0\}$
 is a curve in $\hi.$ Let $C_1^{q_1}+\ldots+C_r^{q_r}$ be its
decomposition into irreducible components. Let $\{P_d=0\}_{red}$ be the reduced curve
$C_1+\ldots+C_r$ and let $\sing(\{P_d=0\}_{red})$ consist of $s$ points $\{ Q_1,\ldots,Q_s\}.$
Suppose that:
\roster
\item  the curve $\{P_{d-k}=0\}$ is reduced;
\item  $Q_i\not\in \{P_{d-k}=0\},\, (i=1,\ldots,s);$ 
\item  for each $j$ with $q_j>1,$ the
curves $C_j$ and $\{P_{d-k}=0\}$ intersect transversally, i.e., the set $C_j\cap\{P_{d-k}=0\}$
consists of $d_j(d-k)$ different points $(d_j=\deg C_j).$
\endroster

The generic fibre of the polynomial $P$ is homotopy equivalent to the bouquet of
$2$-dimensional spheres. In this case the number of these spheres is equal to
$\mu(P)=\dim_\bc \bc[x,y,z]/Jac(P)$ and is equal to 
$$
(d-1)^3-k\cdot\left(\chi(\{P_d=0\})+d(\,2d-{\tilde d}-3)\right)+k^2\cdot(d-{\tilde d}),
$$
where ${\tilde d}=d_1+\ldots+d_r$ is the degree
of the (reduced) curve $\{P_d=0\}_{red},$ \cite{4}.
Let us
consider the following partitioning of the infinite hyperplane $\hi:$ 
\roster
\item the $0$-dimensional stratum $\Xi^0_i$ consisting of one point $Q_i$
each $(i=1,\ldots,s);$
\item the $0$-dimensional stratum $\Lambda_j^0=C_j\cap \{P_{d-k}=0\},$  for each
$j=1,\ldots,r;$
\item  the $1$-dimensional stratum 
$\Xi_j^1=C_j \setminus (\{ Q_i\}\cup \Lambda_j^0),$ for each
$j=1,\ldots,r;$  
\item the $2$-dimensional stratum $\Xi^2=\hi\setminus \{P_d=0\}.$
\endroster
 At
each point of the stratum $\Xi^2,$ the germ of the meromorphic function $P$ has the 
form $\frac{1}{u^d}$ ($\hi=\{u=0\}$). Its infinite zeta-function is equal to $(1-t^d).$ The Euler
characteristic $\chi(\Xi^2)$ of the stratum
$\Xi^2$ is equal to
$$\chi(\hi)-\chi(\{P_d=0\})=3-3{\tilde d}+{\tilde d}^2-
\sum_{i=1}^s\mu_{i},$$
where $\mu_{i}$ is the Milnor
number of the (reduced) curve $\{P_d=0\}_{red}$ at the point $Q_i.$

At each point of the stratum $\Xi_j^1,$ the germ of the meromorphic function $P$ has the 
form $\frac{y_1^{q_j}+u^k}{u^d}.$ 
 Its infinite zeta-function can be determined by the Varchenko type formula from \cite{8} and
is equal to
$$(1-t^{d-k})(1-t^{\frac{q_j(d-k)}{g.c.d.(q_j,k)}})^{-g.c.d.(q_j,k)}.$$
The Euler characteristic of the stratum $\Xi_j^1$ is equal to
$$\chi(C_j)-d_j(d-k)-\sharp\{ C_j\cap \{Q_i:\,i=1,\ldots,s\}\,\}.$$

 At each point of
the stratum $\Lambda_j^0,$ the germ of the meromorphic function $P$ has the 
form $\frac{y_1^{q_j}+u^ky_2}{u^d}.$
Its infinite zeta-function is equal to 1.

 At a point $Q_i,$ the germ of the meromorphic function $P$ has the 
form $\frac{g_i(y_1,y_2)+u^k}{u^d},$ where
$\{g_i=0\}$ is the local equation of the (non-reduced) curve $\{P_d=0\}$ 
at the point $Q_i.$ Its infinite
zeta-function is equal to
$$(1-t^{d-k})\left(\zeta_{g_i}^k(t^{d-k})\right)^{-1}.$$
\remark{Remark 4}
We can not apply the formula $(3)$ directly since the singularity of the germ $g_i$ is, in
general, not isolated. However, it is not difficult to see that, actually, the proof of this formula
uses only the fact that the singularity of the germ $g_i$ can be resolved by a
modification which is an isomorphism outside the origin. This is so for a curve
singularity.
\endremark
 Thus one obtains
$$
\multline
\zeta_P(t) =  (1-t^d)^{\chi(\Xi^2)}
(1-t^{d-k})^{(3{\tilde d}-{\tilde d}^2-{\tilde d}(d-k)+\sum \mu_i)}\,\times \\
\times\prod_{j=1}^r\left(1-t^{\frac{q_j(d-k)}{g.c.d.(q_j,k)}}\right)^{-g.c.d.(q_j,k)\cdot
\chi(\Xi^1_j)}
\cdot\prod_{i=1}^s
\left(\zeta_{g_i}^k(t^{d-k})\right)^{-1}.\qquad\qquad\quad
\endmultline
$$
\bigbreak
\finparrafo

\head\sec On the bifurcation set of a polynomial map
\endhead  

As we have mentioned, a polynomial map $P:\bc^{n+1}\to\bc$ defines a locally
trivial fibration over the complement to a finite set in $\bc.$ The minimal
set $B(P)$ with this property is called the bifurcation set of $P.$ The
bifurcation set consists of critical values of the polynomial $P$ (in the
affine part) and of atypical (``critical") values at infinity. 

 In order to consider a level set $\{P=c\},$ one
can substitute the polynomial $P$ by the polynomial $(P-c)$ and consider
the zero level set. Thus let us consider the zero level set $V_0=\{P=0\}\subset \bc^{n+1}$ of
the polynomial $P.$ 
Let us suppose that the level set $V_0$ of the polynomial $P$ has only isolated
singular  points (in the affine part $\bc^{n+1}).$
For $\rho>0,$ let $B_\rho$ be the open ball of radius
$\rho$ centred at the origin in $\bc^{n+1}$ and  $S_\rho=\partial B_\rho$ be the
$(2n+1)$-dimensional sphere  of radius
$\rho$ with the centre at the origin. There exists $R>0$ such that, for all
$\rho\geq R,$ the sphere $S_\rho$ is transversal to the level set
$V_0=\{P=0\}$ of the polynomial map $P.$  
The restriction $P|_{\bc^{n+1}\setminus B_R}:\bc^{n+1}\setminus B_R\to\bc$ of the function 
$P$ to the complement of the ball $B_R$ defines a $C^\infty$ locally trivial fibration
over a punctured neighbourhood of the origin in $\bc.$
 The loop $\varepsilon_0\cdot\exp(2\pi i
\tau)$ ($0\leq\tau\leq 1,$ $\|\varepsilon_0\|$ small enough) defines the monodromy
transformation $h:V_{\varepsilon_0}\setminus B_R \to
V_{\varepsilon_0}\setminus B_R.$ Let us denote its zeta-function
$\zeta_h(t)$ by $\zeta_P^0(t).$ We use the following definition 

\definition{Definition} The value $0$ is {\it atypical at infinity} for the polynomial $P$
if the restriction
$P|_{\bc^{n+1}\setminus B_R}$ of the function $P$ to the complement of the ball $B_R$ is
not a $C^\infty$ locally trivial fibration over a neighbourhood of the origin in $\bc.$
\enddefinition

\remark{Remark 5} This definition does not depend on a choice of coordinates, i.e.,
it is invariant with respect to polynomial diffeomorphisms of the space $\bc^{n+1}.$
 One can
see that an atypical at infinity value is atypical, i.e. it belongs to the bifurcation set $B(P)$
of the polynomial $P.$ Moreover the bifurcation set $B(P)$ is the union of the set of
critical values of the polynomial $P$ (in $\bc^{n+1}$) and of the set of values atypical at
infinity in the described sense. 
If the singular locus of the level set $V_0=\{P=0\}$ is not finite, the value $0$ hardly can be
considered as typical at infinity. Thus, one should consider this definition as a (possible)
general definition of a value atypical at infinity. In fact the same definition was used in
\cite{10}.
\endremark

Let ${\Cal S}$ be a
prestratification of the infinite hyperplane $\hip$  such that,
for each stratum $\Xi$ of $\Cal S,$ the zero zeta-function
$\zeta_{P,x}^0(t)$ of the germ of the meromorphic function $P$ at a
point $x\in\hip$ does not depend on the point 
$x,$ for
$x
\in
\Xi$ (let it be  $\zeta_{\Xi}^0(t)$ and let its degree be 
$\chi_{\Xi}^0$). 
\proclaim{Theorem 4} 
$$\zeta_{P}^0(t)=\prod_{\Xi\in {\Cal S}}[\zeta_{\Xi}^0(t)]^{\chi
(\Xi)},$$
$$
\chi(V_{\varepsilon_0}\setminus B_R)=\sum_{\Xi\in
{\Cal S}}\chi_{\Xi}^0\cdot\chi(\Xi).$$
\endproclaim

The proof is essentially the same as that of Theorem 1. Since the Euler
characteristic of the set $V_0\setminus B_R$ is equal to $0,$ one has

\proclaim{Corollary  1} If $\zeta_P^0(t)\not\equiv 1,$ then the value $0$ is 
atypical at infinity for the polynomial $P.$
\endproclaim

In several papers (see, e.g.,  \cite{3}, \cite{11}, \cite{12}) there was considered an integer
$\lambda_P(c)\,(c\in\bc)$ such that 
$$\chi(\{P=c\})=\chi(\{P=c+\varepsilon\})+(-1)^{n+1}(\sum
\mu_i+\lambda_P(c) ),$$
\noindent where $\mu_i$  are the Milnor numbers of the (isolated) singular
points of the level set $\{P=c\}\subset\bc^{n+1}.$ Theorem 4 gives the
following formula for this invariant:

\proclaim{Corollary 2}
$$
\lambda_P(0)=(-1)^n\,\deg \zeta_P^0(t)=(-1)^n \sum_{\Xi\in
{\Cal S}}\chi_{\Xi}^0\cdot\chi(\Xi)\,\,\left(=(-1)^n\int_\hip \chi^0_{P,x}\,
d\chi\right).$$
\endproclaim

\example{Example}
Let
$P(x,y,z)=x^ay^b(x^cy^d-z^{c+d})+z,$  $(ad-bc)\not=0,$ and let $\dd=\deg(P)=a+b+c+d.$
The curve $\{P_\dd=0\}\subset\hi$ consists on three components: the line $C_1=\{x=0\}$ with
multiplicity $a,$ the line $C_2=\{y=0\}$ with
multiplicity $b,$ and the reduced curve $C_3=\{x^cy^d-z^{c+d}=0\}.$ Let $Q_1=C_2\cap
C_3=(1:0:0),$
$Q_2=C_1\cap C_3=(0:1:0),$ $Q_3=C_1\cap C_2=(0:0:1).$ At each point $x$ of the infinite
hyperplane $\hi$ except $Q_1$ and $Q_2,$ one has $\zeta_{P,x}^0(t)=1.$ At the point $Q_1,$ the
germ of the meromorphic function $P$ has the form 
$\displaystyle\frac{y^b(y^d-z^{c+d})+zu^{\dd-1}}{u^\dd}.$ Its zero zeta-function can be obtained
by the Varchenko type formula from \cite{8}.  If $(ad-bc)<0,$ then $\zeta_{P,Q_1}^0(t)=1.$
If $(ad-bc)>0,$ then 
$$\zeta_{P,Q_1}^0(t)=(1-t^{\frac{ad-bc}{G.C.D.}})^{G.C.D.},$$
where $G.C.D.=g.c.d(c,d)\cdot g.c.d.(\frac{ad-bc}{g.c.d(c,d)},\dd-1).$ At the point $Q_2,$ we
have just the symmetric situation. Finally 
$$\zeta_{P}^0(t)=(1-t^{\frac{|ad-bc|}{G.C.D.}})^{G.C.D.}.$$
It means that the value $0$ is atypical at infinity. In the same way 
$\zeta_{P-c}^0(t)=1,$ for $c\not=0.$
\endexample

\finparrafo


\enddocument